\documentclass[twoside,11pt]{article}

\usepackage{amsmath,amsfonts}

\newcommand{\D}{\displaystyle}

\begin{document}
\title{{\bf {Optimal decay rate of the bipolar Euler-Poisson system with damping in $\mathbb{R}^3$}}}
\author{Z{\sc higang}
W{\sc u}, \ \  Y{\sc uming} Q{\sc in}\thanks{Corresponding author. E-mail: yuming\_qin@hotmail.com}\\
\footnotesize{\it Department of Applied Mathematics, Donghua
University, Shanghai, 201620, P.R. China}}

\date{}
\maketitle

\textbf{{\bf Abstract:}} By rewriting a bipolar Euler-Poisson
equations with damping into an Euler equation with damping coupled
 with an Euler-Poisson equation with damping, and using a new
spectral analysis, we obtain the optimal decay results of the
solutions in $L^2$-norm, which improve theose in \cite{Li3, Wu3}.
More precisely, the velocities $u_1,u_2$ decay at the $L^2-$rate
$(1+t)^{-\frac{5}{4}}$, which is faster than the normal $L^2-$rate
$(1+t)^{-\frac{3}{4}}$ for the Heat equation and the Navier-Stokes
equations. In addition, the disparity of two densities
$\rho_1-\rho_2$ and the disparity of two velocities $u_1-u_2$ decay
at the $L^2$-rate $(1+t)^{-2}$.

 \bigbreak \textbf{{\bf Key Words}:}
Bipolar damped Euler-Poisson system; decay estimates.

\bigbreak
\section*{ 1.\ Introduction}

\quad\quad The compressible bipolar Euler-Poisson equations with
damping (BEP) takes the following form
$$
\left\{\begin{array}{l}
\partial_t\rho_1+{\rm div} (\rho_1 u_1)=0, \\
\partial_t(\rho_1 u_1)+ {\rm div}(\rho_1u_1\otimes u_1)+\nabla
P(\rho_1)
=\rho_1\nabla\phi-\rho_1u_1, \\
\partial_t\rho_2+{\rm div} (\rho_2 u_2)=0, \\
\partial_t(\rho_2 u_2)+ {\rm div}(\rho_2u_2\otimes u_2)+\nabla
P(\rho_2)
=-\rho_2\nabla\phi-\rho_2u_2, \\
\Delta\phi=\rho_1-\rho_2,\ \ \ x\in\mathbb{R}^3,\ t\geq0.
\end{array}
        \right.
        \eqno({\rm 1.1})
$$
Here the unknown functions $\rho_i(x,t),u_i(x,t)\ (i=1,2),\phi(x,t)$
are the charge densities, current densities, velocities and
electrostatic potential, respectively, and the pressure
$P=P(\rho_i)$ is a smooth function with $P'(\rho_i)>0$ for
$\rho_i>0$. The system (1.1) is described charged particle fluids,
for example, electrons and holes in semiconductor
 devices, positively and negatively charged ions in a plasma. We
 refer to \cite{Jugel,Markowich,Sitenko} for the physical background
 of the system (1.1).

In this paper, we want to study the system (1.1) with the Cauchy
data
$$
\rho_i(x,0)=\rho_{i0}(x)>0,\ u_i(x,0)=u_{i0}(x),\ i=1,2. \eqno(1.2)
$$

A lot of important works for the system (1.1) have been done. For
one-dimensional case, we refer to Zhou and Li \cite{Zhou}, and Tsuge
\cite{Tsuge} for the unique existence of the stationary solutions,
Natalini \cite{Natalini1} and Hsiao and Zhang \cite{Hsiao1} for
global entropy weak solutions on the whole real line and bounded
domain respectively,  Huang and Li \cite{Huang} for the large-time
behavior and quasi-neutral limit of $L^\infty$-solution, Natalini
\cite{Natalini1}, Hsiao and Zhang \cite{Hsiao2} for the
relaxation-time limit, Gasser and Marcati \cite{Gasser2} for the
combined limit, Zhu and Hattori \cite{Zhu} for the stability of
steady-state solutions to a recombined one-dimensional bipolar
hydrodynamical model, Gasser, Hsiao and Li \cite{Gasser1} for
large-time behavior of smooth solutions with small initial data.

For the multi-dimensional case, Lattanzio \cite{Lattanzio} studied
the relaxation limit, and Li \cite{Li2} investigated the diffusive
relaxation. Ali and J\"{u}ngel \cite{Ali} and Li and Zhang
\cite{Li1} considered the global smooth solutions of the Cauchy
problem in the Sobolev's space and Besov's space, respectively.
Later, Ju \cite{Ju} investigated the global existence of smooth
solution to the IBVP for the system (1.1).

Recently, by using the standard energy method together with the
analysis of the Green's function, Li and Yang \cite{Li3}
investigated the decay rate of the Cauchy's problem of the system
(1.1) of the classical solution when the initial data are small in
the space $H^3\cap L^1$. They deduced that the electric field (a
nonlocal term in hyperbolic-parabolic system) slows down the decay
rate of the velocity of the BEP system. More precisely, they
obtained $\|\rho_i-\bar{\rho}\|_{L^2}\sim(1+t)^{-\frac{3}{4}},\
\|u_i\|_{L^2}\sim(1+t)^{-\frac{1}{4}}$. We refer to the relevant
works on the unipolar Navier-Stokes-Poisson equations (NSP) and
unipolar Euler-Poisson equations with damping
\cite{LiH1,LiH2,Zhang,Wang1,Wu1,Wu2}. The decay result in \cite{Li3}
was improved by us in \cite{Wu3} as
$\|u_i\|_{L^2}\sim(1+t)^{-\frac{3}{4}}$ by using the method
introduced by Guo and Wang \cite{Guo,WangY1}, which is based on a
family of energy estimates with minimum derivative counts and
interpolations among them without linear decay analysis. The
additional assumption on the initial perturbation in \cite{Wu3} is
$\nabla\phi_0\in L^2\ (\rho_0-\bar{\rho}\in \dot{H}^{-1})$, but it
need not the smallness of the initial perturbation in $L^1$-norm.

On the other hand, recently, Tan and Wu \cite{Tan} gave a new and
detailed spectral analysis for the Euler equations with damping by
taking Hodge decomposition on the velocity. In fact, they deduced
the decay rate: $\|\rho-\bar{\rho}\|_{L^2}\sim(1+t)^{-\frac{3}{4}},\
\|u\|_{L^2}\sim(1+t)^{-\frac{5}{4}}$, which improved the former
decay rate in Wang and Yang \cite{Wang}.

The purpose of this paper is to improve the $L^2$-norm decay
estimates in Li and Yang \cite{Li3} and Wu and Wang \cite{Wu3}.
Comparing with \cite{Li3,Wu3}, our essential idea is to change the
system (1.1) into an Euler equation with damping coupled with an
Euler-Poisson equation with damping. Our another observation is that
when two densities $\rho_1$ and $\rho_2$ are the perturbations of
$\bar{\rho}$ and $u_1$ and $u_2$ are the perturbations of zero, then
$\rho_1-\rho_2$ and $u_1-u_2$ can be small in $H^3\cap L^1$. Hence,
by the spectral analysis in \cite{Tan} for the Euler equations with
damping and the spectral analysis in Wu and Wang \cite{Wu2} for the
Euler-Poisson equations with damping, together with the global
existence and energy estimates in \cite{Wu3}, we can achieve
expected decay results for the rewritten system (2.1). As a
byproduct, we obtain the improved decay results of the solution for
the original system (1.1).

Our main results are stated in the following theorems.

\bigbreak \noindent\textbf{Theorem 1.1.}(\cite{Wu3}) Let
$P'(\rho_i)>0 (i=1,2)$ for $\rho_i>0$, and $\bar{\rho}>0$. Assume
that $(\rho_i-\bar{\rho},u_{i0})\in H^3(\mathbb{R}^3)$ for $i=1,2$,
with
$\epsilon_0=:\|(\rho_{i0}-\bar{\rho},u_{i0})\|_{H^3(\mathbb{R}^3)}$
small. Then there is a unique global classical solution
$(\rho_i,u_i,\nabla\phi)$ of the Cauchy problem (1.1)-(1.2)
satisfying
$$
\|(\rho_i-\bar{\rho},u_i,\nabla\phi)\|_{H^3}^2\leq
C\epsilon_0.\eqno(1.3)
$$

\bigbreak \noindent\textbf{Theorem 1.2.} Under the assumptions of
Theorem 1.1. If further,
$$
\|(\rho_{10}-\rho_{20},\rho_{10}+\rho_{20}-2\bar{\rho},u_{10},u_{20})\|_{L^1}\leq\epsilon_1\
(\epsilon_1\ll1),\eqno(1.4)
$$
we have
$$
\|\partial_x^k(\rho_1+\rho_2-2\bar{\rho})\|_{L^2}\leq
C(\epsilon_0+\epsilon_1)(1+t)^{-\frac{3}{4}-\frac{|k|}{2}},\
k=0,1,\eqno(1.5)
$$
$$
\|\partial_x^k(u_1+u_2)\|_{L^2}\leq
C(\epsilon_0+\epsilon_1)(1+t)^{-\frac{5}{4}-\frac{|k|}{2}},\
k=0,1,\eqno(1.6)
$$
$$
\|(\rho_1-\rho_2,u_1-u_2)\|_{L^2}\leq
C(\epsilon_0+\epsilon_1)(1+t)^{-2},\eqno(1.7)
$$
$$
\|\partial_x(\rho_1-\rho_2,u_1-u_2)\|_{L^2}\leq
C(\epsilon_0+\epsilon_1)(1+t)^{-\frac{15}{8}}.\eqno(1.8)
$$

\bigbreak \noindent\textbf{Corollary 1.1.} Under the assumptions of
Theorem 1.2, we have
$$
\|\partial_x^k(\rho_{1}-\bar{\rho},\rho_{2}-\bar{\rho})\|_{L^2}\leq(\epsilon_0+\epsilon_1)(1+t)^{-\frac{3}{4}-\frac{|k|}{2}},\
k=0,1,\eqno(1.9)
$$
and
$$
\|\partial_x^k(u_1,u_2)\|_{L^2}\leq(\epsilon_0+\epsilon_1)(1+t)^{-\frac{5}{4}-\frac{|k|}{2}},\
k=0,1.\eqno(1.10)
$$

\bigbreak \noindent\textbf{Remark 1.1.} The decay rates in Corollary
1.1 are optimal. This can be understood from Theorem 1.2 in Li et\
al. \cite{LiH0}, where the same decay rates as (1.9) and (1.10) from
below were shown.

\bigbreak \noindent\textbf{Remark 1.2.} From Corollary 1.1, we know
the density $\rho_i-\bar{\rho}$ have the same decay rate in
$L^2$-norm as the solution of the Navier-Stokes equations, while the
velocity $u_i$ in \cite{Li3} decays faster than the densities. That
is, we have improved the decay result in \cite{Li3,Wu3}. In
addition, the decay rates of the disparity of two densities
$\rho_1-\rho_2$ and the disparity of two velocities $u_1-u_2$ in
(1.7) are surprising and satisfactory, because of the exponential
decay rate for the Euler-Poisson equations with damping and the
coupling of two systems in (2.4).

\bigbreak \noindent\textbf{Remark 1.3.}   The authors in
\cite{Duan1,Wu3} obtained the optimal decay rate based on some
Lyapunov energy functionals and a family of energy estimates with
minimum derivative counts, respectively. Consequently, both
\cite{Duan1} and \cite{Wu3} need not the smallness of $L^1$-norm of
the initial data. However, the usual energy method can not be used
to the rewritten system (2.4), thus, it is also interesting to get
rid of the smallness of the assumption
$\|(\rho_{10}-\rho_{20},\rho_{10}+\rho_{20}-2\bar{\rho},u_{10},u_{20})\|_{L^1}$
in Theorem 1.2.

\bigbreak \noindent\textbf{Remark 1.4.} The similar problem for the
bipolar quantum hydrodynamic model will be investigated in a
forthcoming paper.

\bigbreak \noindent\textbf{Notations.} In this paper,
$D^l=\partial_x^l$ with an integer $l\geq 0$ stands for the usual
any spatial derivatives of order $l$. For $1\leq p\leq \infty$ and
an integer $m\geq 0$, we use $L^p$ and $W^{m,p}$ denote the usual
Lebesgue space $L^p(\mathbb{R}^3)$ and Sobolev spaces
$W^{m,p}(\mathbb{R}^3)$ with norms $\|\cdot\|_{L^p}$ and
$\|\cdot\|_{W^{m,p}}$, respectively, and set $H^m=W^{m,2}$ with norm
$\|\cdot\|_{H^m}$ when $p=2$. In addition, we define the homogeneous
Sobolev's space $\dot{H}^s$ of all $g$ for which $\|g\|_{\dot{H}^s}$
is finite, where
$$
\|g\|_{\dot{H}^s}:=\|\Lambda^sg\|_{L^2}=\||\xi|^s\hat{g}\|_{L^2}.
$$
Here, for $s\in\mathbb{R}$, a pseudo-differential operator
$\Lambda^s$ by
$$
\Lambda^sg(x)=\int_{\mathbb{R}^n}|\xi|^s\hat{g}(\xi){\rm
e}^{2\pi\sqrt{-1}x\cdot\xi} {\rm d}\xi,
$$
where $\hat{g}$ denotes the Fourier transform of $g$.

Throughout this paper, for simplicity, we write
$\|\cdot\|_{L^\infty}=\|\cdot\|_\infty,\ \|\cdot\|_{L^2}=\|\cdot\|$.
And we will use $C$ or $C_i$ denotes a positive generic (generally
large) constant that may vary at different places. And the sign
``$\sim$" means
$$
f\sim g \Longleftrightarrow {\rm there\ exists\ two\ positive\
constants}\ C_1, C_2\ {\rm such\ that\ }\ C_1|f|\leq |g|\leq C_2|f|.
$$

The rest of the paper is arranged as follows. In Section 2, we
reformulate the original system, then give detailed spectral
analysis for the Euler equations with damping and Euler-Poisson
equations with damping. The proof of temporal decay results of the
solution will be derived in Section 3.

\section*{ 2.\ Spectral analysis and linear $L^2$ estimates}

\subsection*{ 2.1\ Reformulation of original problem}

\quad\quad Assume $\bar{\rho}=1$ and $P'(\bar{\rho})=1$ without loss
of generality. Then the Cauchy problem (1.1)-(1.2) is reformulated
as
$$
\left\{\begin{array}{l}
\partial_t \rho_1+{\rm div} u_1=-u_1\cdot\nabla \rho_1-\rho_1{\rm div}u_1, \\[1mm]
\partial_t u_1+u_1+\nabla \rho_1-\nabla\phi=-u_1\cdot\nabla
u_1-h(\rho_1)\nabla \rho_1, \\[1mm]
\partial_t \rho_2+{\rm div} u_2=-u_2\cdot\nabla \rho_2-\rho_2{\rm div}u_2, \\[1mm]
\partial_t u_2+u_2+\nabla \rho_2+\nabla\phi=-u_2\cdot\nabla
u_2-h(\rho_2)\nabla \rho_2, \\[1mm]
\Delta\phi=\rho_1-\rho_2,\\[1mm]
(\rho_1,u_1,\rho_2,u_2)(x,0)=(\rho_{10},u_{10},\rho_{20},u_{20})(x),
\end{array}
        \right.
        \eqno({\rm 2.1})
$$
where $h(\rho_i)=\frac{P'(\rho_i)}{\rho_i}-1$.

Next, set
$$
n_1=\rho_1+\rho_2-2,\ n_2=\rho_1-\rho_2,\ w_1=u_1+u_2,\
w_2=u_1-u_2,\eqno(2.2)
$$
which give equivalently
$$
\rho_1=\frac{n_1+n_2}{2}+1,\ \rho_2=\frac{n_1-n_2}{2}+1,\
u_1=\frac{w_1+w_2}{2},\ u_2=\frac{w_1-w_2}{2}.\eqno(2.3)
$$

From (2.2) and (2.3), it follows that the Cauchy problem (2.1) can
be reformulated into the following Cauchy problem for the unknown
$(n_1,w_1,n_2,w_2,\phi)$
$$
\left\{\begin{array}{l}
\partial_t n_1+{\rm div} w_1=f_1(n_1,w_1,n_2,w_2), \\[1mm]
\partial_t w_1+w_1+\nabla n_1=f_2(n_1,w_1,n_2,w_2), \\[1mm]
\partial_t n_2+{\rm div} w_2=f_3(n_1,w_1,n_2,w_2), \\[1mm]
\partial_t w_2+w_2+\nabla n_2-2\nabla\phi=f_4(n_1,w_1,n_2,w_2), \\[1mm]
\Delta\phi=n_2,\\[1mm]
(n_1,w_1,n_2,w_2)(x,0)=(n_{10},w_{10},n_{20},w_{20})(x),
\end{array}
        \right.
        \eqno({\rm 2.4})
$$
where
$(n_{10},w_{10},n_{20},w_{20}):=(\rho_{10}+\rho_{20}-2,u_{10}+u_{20},\rho_{10}-\rho_{20},u_{10}-u_{20})$,
and
$$  \begin{array}[b]{rl}
 f_1=&\D-\frac{1}{2}[w_2\nabla n_1+w_1\nabla n_2+n_1{\rm
div}w_1+n_2{\rm div}w_2];\\[2mm]
 f_2=&\D-\frac{1}{2}[w_1\nabla w_1+w_2\nabla
w_2+(h(\frac{n_1\!+\!n_2}{2}+1)\!+\!h(\frac{n_1\!-\!n_2}{2}+1))\nabla
n_1\\[2mm]
&\D+(h(\frac{n_1\!+\!n_2}{2}+1)-h(\frac{n_1\!-\!n_2}{2}+1))\nabla n_2];\\[2mm]
 f_3=&\D-\frac{1}{2}[w_1\nabla n_2\!+\!w_2\nabla n_1+n_1{\rm
div}w_2+n_2{\rm div}w_1];\\[2mm]
 f_4=&\D-\frac{1}{2}[w_1\nabla w_2+w_2\nabla
w_1+(h(\frac{n_1\!+\!n_2}{2}+1)\!-\!h(\frac{n_1\!-\!n_2}{2}+1))\nabla
n_1\\[2mm]
&\D+(h(\frac{n_1\!+\!n_2}{2}+1)\!+\!h(\frac{n_1\!-\!n_2}{2}+1))\nabla
n_2].
\end{array}
\eqno(2.5)
$$

From the foregoing, the Cauchy problem (2.4) can be formally divided
into the Cauchy problem for the Euler equations with damping
$(2.4)_{1,2}$ and the Euler-Poisson equations with damping
$(2.4)_{3,4,5}$, which interact each other through the nonlinear
inhomogeneous terms on the right-hand side.

\subsection*{ 2.2\ Spectral analysis of linearized Euler equations with damping}

\quad\quad In this subsection, we shall give a spectral analysis for
the linearized Euler equations with damping. A similar analysis can
be founded in \cite{Tan}. For the convenience of the readers, the
proof will be also sketched here.

We take Hodge decomposition to analyze the following linearized
equations corresponding to system $(2.4)_{1,2}$
$$
\left\{\begin{array}{l}
\partial_t n_1+{\rm div} w_1=0, \\[1mm]
\partial_t w_1+w_1+\nabla n_1=0,\\[1mm]
(n_1,w_1)(t=0)=(n_{10},w_{10}).
\end{array}
        \right.
        \eqno({\rm 2.6})
$$

Let $v_1=\Lambda^{-1}{\rm div}w_1$ be the ``compressible part" of
the velocity and $d_1=\Lambda^{-1}{\rm curl}w_1$ be the
``incompressible part", then the system (2.6) can be rewritten as
$$
\left\{\begin{array}{l}
\partial_t n_1+\Lambda v_1=0, \\[1mm]
\partial_t v_1+w_1+\Lambda n_1=0,\\[1mm]
\partial_t d_1+d_1=0.
\end{array}
        \right.
        \eqno({\rm 2.7})
$$

Obviously, by the definition of $v_1$ and $d_1$, and the relation
$$
w_1=-\Lambda^{-1}\nabla v_1-\Lambda^{-1}{\rm div}d_1,\eqno(2.8)
$$
the estimates in space $H^l(\mathbb{R}^3)$ for the original function
$w_1$ is the same as for the functions $(v_1,d_1)$.

By the semigroup theory for evolutionary equation, the solution
$(n_1,v_1)$ of the first two equations of the system (2.7) can be
expressed via the Cauchy problem for $X=(n_1,v_1)^\tau$ as
$$
X_t=AX,\ X(0)=X_0,\ t\geq0,
$$
which leads to
$$
X(t)=S_1(t)X_0=:e^{tA}X_0,\ t\geq0.
$$

Taking the Fourier transform with respect to the space variable
yields
$$
\frac{d}{dt}\hat{X}=A(\xi)\hat{X}\ {\rm with}\ A(\xi)=\left(
                                                        \begin{array}{cc}
                                                          0 & -|\xi| \\
                                                          |\xi| & -1 \\
                                                        \end{array}
                                                      \right).
\eqno(2.9)
$$

The characteristic polynomial of $A(\xi)$ is
$\lambda^2+\lambda+|\xi|^2$ and has two distinct roots:
$$
\lambda_\pm(\xi)=-\frac{1}{2}\pm\frac{1}{2}\sqrt{1-|\xi|^2}.\eqno(2.10)
$$

By a direct computation, we can verify the exact expression the
Fourier transform $\hat{G_1}(\xi,t)$ of Green function
$G_1(x,t)=e^{tA}$ as
$$
\hat{G_1}(\xi,t)=\left(
                 \begin{array}{cc}
                   \frac{\lambda_+e^{\lambda_-t}-\lambda_-e^{\lambda_+t}}{\lambda_+-\lambda_-} & |\xi|\frac{e^{\lambda_+t}-e^{\lambda_-t}}{\lambda_+-\lambda_-} \\
                   -|\xi|\frac{e^{\lambda_+t}-e^{\lambda_-t}}{\lambda_+-\lambda_-} & \frac{\lambda_+e^{\lambda_-t}-\lambda_-e^{\lambda_+t}}{\lambda_+-\lambda_-}-\frac{e^{\lambda_+t}-e^{\lambda_-t}}{\lambda_+-\lambda_-} \\
                 \end{array}
               \right).
               \eqno(2.11)
               $$

We use the standard higher-lower frequency decomposition to derive
the long-time decay rate of solutions in $L^2$ framework. For
$|\xi|\ll1$, it holds that
$$\begin{array}[b]{rl}
\frac{\lambda_+e^{\lambda_-t}-\lambda_-e^{\lambda_+t}}{\lambda_+-\lambda_-}\sim
e^{-|\xi|^2t}-|\xi|^2e^{-t},\\[1mm]
|\xi|\frac{e^{\lambda_+t}-e^{\lambda_-t}}{\lambda_+-\lambda_-}\sim
-|\xi|e^{-t}+|\xi|e^{-|\xi|^2t},\\[1mm]
\frac{\lambda_+e^{\lambda_-t}-\lambda_-e^{\lambda_+t}}{\lambda_+-\lambda_-}-\frac{e^{\lambda_+t}-e^{\lambda_-t}}{\lambda_+-\lambda_-}
\sim (1-|\xi|^2)e^{-t}. \end{array} \eqno(2.12)
$$

On the other hand, as we know, the higher frequency part of
$\hat{G}_1(\xi,t)$ in $L^2$-norm has the exponential decay rate.

Then by using the explicit expression of $\hat{n_1}$ and
$\hat{v_1}$:
$$
\hat{n_1}=\frac{\lambda_+e^{\lambda_-t}-\lambda_-e^{\lambda_+t}}{\lambda_+-\lambda_-}\hat{n_{10}}
+|\xi|\frac{e^{\lambda_+t}-e^{\lambda_-t}}{\lambda_+-\lambda_-}\hat{v_{10}},
$$
$$
\hat{v_1}=-|\xi|\frac{e^{\lambda_+t}-e^{\lambda_-t}}{\lambda_+-\lambda_-}\hat{n_{10}}
+\frac{\lambda_+e^{\lambda_-t}-\lambda_-e^{\lambda_+t}}{\lambda_+-\lambda_-}-\frac{e^{\lambda_+t}-e^{\lambda_-t}}{\lambda_+-\lambda_-}\hat{v_{10}},
$$
we have
$$\begin{array}[b]{rl}
\|\hat{n_1}(t)\|^2=&\D\int_{|\xi|\leq\eta}|\hat{n_1}(\xi,t)|^2d\xi+\int_{|\xi|\geq\eta}|\hat{n_1}(\xi,t)|^2d\xi\\[1mm]
\leq &\D
C\int_{|\xi|\leq\eta}e^{-|\xi|^2t}(|\hat{n_{10}}|^2+|\hat{v_{10}}|^2)d\xi+Ce^{-bt}\int_{|\xi|\geq\eta}(|\hat{n_{10}}|^2+|\hat{v_{10}}|^2)d\xi\\[1mm]
\leq &\D
C\|(n_{10},v_{10})\|_{L^1}^2\int_{|\xi|\leq\eta}e^{-|\xi|^2t}d\xi+Ce^{-bt}\|(n_{10},v_{10})\|_{L^2}^2\\[1mm]
\leq &\D C(1+t)^{-\frac{3}{2}}\|(n_{10},v_{10})\|_{L^2\cap L^1}^2,\
{\rm where}\ b>0,
\end{array}\eqno(2.13)
$$
and
$$\begin{array}[b]{rl}
\|\hat{v_1}(t)\|^2=&\D\int_{|\xi|\leq\eta}|\hat{v_1}(\xi,t)|^2d\xi+\int_{|\xi|\geq\eta}|\hat{v_1}(\xi,t)|^2d\xi\\[1mm]
\leq &\D
C\int_{|\xi|\leq\eta}|\xi|^2e^{-|\xi|^2t}(|\hat{n_{10}}|^2+|\hat{v_{10}}|^2)d\xi+Ce^{-bt}\int_{|\xi|\geq\eta}(|\hat{n_{10}}|^2+|\hat{v_{10}}|^2)d\xi\\[1mm]
\leq &\D
C\|(n_{10},v_{10})\|_{L^1}^2\int_{|\xi|\leq\eta}|\xi|^2e^{-|\xi|^2t}d\xi+Ce^{-bt}\|(n_{10},v_{10})\|_{L^2}^2\\[1mm]
\leq &\D C(1+t)^{-\frac{5}{2}}\|(n_{10},v_{10})\|_{L^2\cap L^1}^2,\
{\rm where}\ b>0.
\end{array}
\eqno(2.14)
$$

Similarly, the $L^2$-decay rate for the derivatives of $\hat{n_1}$
and $\hat{v_1}$ are
$$
\|\partial_x^kn_1\|_{L^2}\leq
C(1+t)^{-\frac{3}{4}-\frac{k}{2}}\|n_{10},v_{10}\|_{\dot{H}^k\cap
L^1},\eqno(2.15)
$$
and
$$
\|\partial_x^kv_1\|_{L^2}\leq
C(1+t)^{-\frac{5}{4}-\frac{k}{2}}\|n_{10},v_{10}\|_{\dot{H}^k\cap
L^1}.\eqno(2.16)
$$

From (2.15) and (2.16), and using the relation $
w_1=-\Lambda^{-1}\nabla v_1-\Lambda^{-1}{\rm div}d_1$ and the fact
$\|d_1\|_{L^2}\sim e^{-t}$, one can easily obtain the following
$L^2$-decay result for the linearized Euler equation with damping.

\bigbreak \noindent\textbf{Lemma 2.1.} Let $(n_{10},w_{10})\in
H^l\cap L^1$ and $n_1,w_1$ satisfy the system (2.6). Then there
exists a constant $C$ such that for $0\leq k\leq l$
$$
\|\partial_x^kn_1\|_{L^2}\leq
C(1+t)^{-\frac{3}{4}-\frac{k}{2}}\|(n_{10},w_{10})\|_{\dot{H}^k\cap
L^1},\eqno(2.17)
$$
and
$$
\|\partial_x^kw_1\|_{L^2}\leq
C(1+t)^{-\frac{5}{4}-\frac{k}{2}}\|(n_{10},w_{10})\|_{\dot{H}^k\cap
L^1}.\eqno(2.18)
$$

\subsection*{ 2.3\ Spectral analysis of linearized Euler-Poisson equations with damping}

\quad\quad In this subsection, we shall give a spectral analysis for
the linearized Euler-Poisson equations with damping.

We take Hodge decomposition to analyze the following linearized
equations corresponding to system $(2.4)_{3,4,5}$
$$
\left\{\begin{array}{l}
\partial_t n_2+{\rm div} w_2=0, \\[1mm]
\partial_t w_2+w_2+\nabla n_2-2\nabla \phi=0,\\[1mm]
\Delta\phi=n_2,\\[1mm]
(n_2,w_2)(t=0)=(n_{20},w_{20}).
\end{array}
        \right.
        \eqno({\rm 2.19})
$$

Let $v_2=\Lambda^{-1}{\rm div}w_2$ be the ``compressible part" of
the velocity and $d_2=\Lambda^{-1}{\rm curl}w_2$ be the
``incompressible part", then the system (2.19) can be rewritten as
$$
\left\{\begin{array}{l}
\partial_t n_2+\Lambda v_2=0, \\[1mm]
\partial_t v_2+w_2-\Lambda n_2+2\Lambda^{-1} n_2=0,\\[1mm]
\partial_t d_2+d_2=0.
\end{array}
        \right.
        \eqno({\rm 2.20})
$$

By the semigroup theory for evolutionary equations, the solution
$(n_2,v_2)$ of the first two equations of the system (2.20) can be
expressed via the Cauchy problem for $Y=(n_2,v_2)^\tau$ as
$$
Y_t=BY,\ Y(0)=Y_0,\ t\geq0,
$$
which leads to
$$
Y(t)=S_1(t)Y_0=:e^{tB}Y_0,\ t\geq0.
$$

Taking the Fourier transform with respect to the space variable
yields
$$
\frac{d}{dt}\hat{Y}=B(\xi)\hat{Y}\ {\rm with}\      B(\xi)=\left(
                                                        \begin{array}{cc}
                                                          0 & -|\xi| \\
                                                          |\xi|-2|\xi|^{-1} & -1 \\
                                                        \end{array}
                                                      \right).
\eqno(2.21)
$$

The characteristic polynomial of $B(\xi)$ is
$\lambda^2+\lambda+|\xi|^2+2$ and has two distinct roots:
$$
\lambda_\pm(\xi)=-\frac{1}{2}\pm\frac{1}{2}\sqrt{-7-4|\xi|^2}.\eqno(2.22)
$$

By a direct computation, we can verify the exact expression the
Fourier transform $\hat{G_2}(\xi,t)$ of Green function
$G_2(x,t)=e^{tB}$ as
$$
\hat{G_2}(\xi,t)=\left(
                 \begin{array}{cc}
                   \frac{\lambda_+e^{\lambda_-t}-\lambda_-e^{\lambda_+t}}{\lambda_+-\lambda_-} & |\xi|\frac{e^{\lambda_+t}-e^{\lambda_-t}}{\lambda_+-\lambda_-} \\
                   -|\xi|\frac{e^{\lambda_+t}-e^{\lambda_-t}}{\lambda_+-\lambda_-} & \frac{\lambda_+e^{\lambda_-t}-\lambda_-e^{\lambda_+t}}{\lambda_+-\lambda_-}-\frac{e^{\lambda_+t}-e^{\lambda_-t}}{\lambda_+-\lambda_-} \\
                 \end{array}
               \right).
               \eqno(2.23)
               $$
We use the standard higher-lower frequency decomposition to derive
the long-time decay rate of solutions in $L^2$ framework. For
$|\xi|\ll1$, it holds that
$$
\begin{array}[b]{rl}
\frac{\lambda_+e^{\lambda_-t}-\lambda_-e^{\lambda_+t}}{\lambda_+-\lambda_-}\sim
-(1+|\xi|^2)e^{-t},\\[1mm]
|\xi|\frac{e^{\lambda_+t}-e^{\lambda_-t}}{\lambda_+-\lambda_-}\sim
-|\xi|e^{-t},\\[1mm]
\frac{\lambda_+e^{\lambda_-t}-\lambda_-e^{\lambda_+t}}{\lambda_+-\lambda_-}-\frac{e^{\lambda_+t}-e^{\lambda_-t}}{\lambda_+-\lambda_-}
\sim e^{-t}. \end{array} \eqno(2.24)
$$

On the other hand, it is obvious that the higher frequency part of
$\hat{G}(\xi,t)$ in $L^2$-norm has also the exponential decay rate.

Then by using the explicit expression of $\hat{n_2}$ and
$\hat{v_2}$, we have the following $L^2$-decay rate for $\hat{n_2}$
and $\hat{v_2}$:
$$
\|\partial_x^k(n_2,w_2)\|_{L^2}\leq
Ce^{-bt}\|(n_{10},v_{10})\|_{\dot{H}^k\cap L^1},\ {\rm for\ some\
constant}\ b>0.\eqno(2.25)
$$
From (2.25) and the relation $ w_2=-\Lambda^{-1}\nabla
v_2-\Lambda^{-1}{\rm div}d_2$ and the fact $\|d_2\|_{L^2}\sim
e^{-t}$, one can easily obtain the following $L^2$-decay result for
the linearized Euler equation with damping.

\bigbreak \noindent\textbf{Lemma 2.2.} Let $(n_{20},w_{20})\in
H^l\cap L^1$ and $n_2,w_2$ satisfy the system (2.19). Then there
exists a positive constant $C$ such that for $0\leq k\leq l$
$$
\|\partial_x^k(n_2,w_2)\|_{L^2}\leq
Ce^{-bt}\|(n_{20},w_{20})\|_{\dot{H}^k\cap L^1},\ {\rm for\ some\
constant}\ b>0. \eqno(2.26)
$$

\section*{ 3.\ Proof of the main results}

\quad\quad In this section, we will mainly use Lemma 2.1 and Lemma
2.2 to prove the optimal decay rate of the nonlinear system (2.4).

Recall that
$$
n_1=\rho_1+\rho_2-2\bar{\rho},\ n_2=\rho_1-\rho_2,\ w_1=u_1+u_2,\
w_2=u_1-u_2.\eqno(3.1)
$$

Then the assumption (1.4) implies that
$$
\|(n_{10},n_{20},w_{10},w_{20})\|_{L^1}\leq \epsilon_1,\
\epsilon_1\ll1.\eqno(3.2)
$$

On the other hand, the global existence in Theorem 1.1 implies that
$$
\frac{1}{2}\leq n_1,n_2\leq 2,\ |h(\frac{n_1\pm n_2}{2}+1)|\leq
C|n_1\pm n_2|,\ |h^{(k)}(\cdot)|\leq C,\ i=1,2.\eqno(3.3)
$$

Now, we are in a position to prove Theorem 1.2. Define
$$
 \begin{array}[b]{rl}
M(t)\!=\! \D\!\!&\sup\limits_{0\leq s\leq t}
\{\sum\limits_{k=0,1}[(1\!+\!s)^{\frac{3}{4}\!+\!\frac{k}{2}}\|D^kn_1\|_{L^2}\!+\!(1\!+\!s)^{\frac{5}{4}\!+\!\frac{k}{2}}\|D^kw_1\|_{L^2}
\!+\!(1\!+\!s)^2\|(n_2,w_2)\|_{L^2}]\\
 &\D+(1\!+\!s)^{\frac{15}{8}}\|D(n_2,w_2)\|_{L^2}+(1\!+\!s)^{\frac{5}{4}}\|D^2(n_1,w_1,n_2,w_2)\|_{L^2}+\|D^3(n_1,w_1,n_2,w_2)\|_{L^2}\}.
 \end{array}
 \eqno(3.4)
 $$
We claim that it holds for any $t\in[0,T]$,
$$
M(t)\leq C(\epsilon_1+\epsilon_2).\eqno(3.5)
$$

First, by Duhamel principle and Lemma 2.1, Lemma 2.2, it is easy to
verity that $L^2$-norm of the solution $n_1,n_2,w_1,w_2$ of the
problem (2.1) can be expressed as
$$
 \begin{array}[b]{ll}
\|D^kn_1\|(t)\leq
\!C\|n_{10},w_{10}\|_{L^1}(1\!+\!t)^{-(\frac{3}{4}\!+\!\frac{k}{2})}\!+\!\int_0^t(1\!+\!t\!-\!s)^{-(\frac{3}{4}\!+\!\frac{k}{2})}[\|(f_1,f_2)\|_{L^1}\!+\!\|D^{k}(f_1,f_2)\|](s)ds,\\[1mm]
\|D^kw_1\|(t)\leq
\!C\|n_{10},w_{10}\|_{L^1}(1\!+\!t)^{-(\frac{5}{4}\!+\!\frac{k}{2})}\!+\!\int_0^t(1\!+\!t\!-\!s)^{-(\frac{5}{4}\!+\!\frac{k}{2})}[\|(f_1,f_2)\|_{L^1}\!+\!\|D^{k}(f_1,f_2)\|](s)ds,\\[1mm]
\|D^kn_2\|(t)\leq
\!Ce^{-bt}\|n_{20},w_{20}\|_{L^1}+\int_0^te^{-b(t-s)}[\|(f_3,f_4)\|_{L^1}\!+\!\|D^{k}(f_3,f_4)\|](s)ds,\\[1mm]
\|D^kw_2\|(t)\leq
\!Ce^{-bt}\|n_{20},w_{20}\|_{L^1}+\int_0^te^{-b(t-s)}[\|(f_3,f_4)\|_{L^1}\!+\!\|D^{k}(f_3,f_4)\|](s)ds,\
b>0.
\end{array}
\eqno(3.6)
$$

\bigbreak \noindent\textbf{Step 1. Basic estimates} From (3.3), the
nonlinear terms $f_1,f_2,f_3,f_4$ can be rewritten into
$$
 \begin{array}[b]{rl}
f_1\sim &\D\mathcal{O}(1)(w_2Dn_1+w_1Dn_2+n_1Dw_1+n_2Dw_2),\\[1mm]
f_2\sim &\D
\mathcal{O}(1)(w_1Dw_1+w_2Dw_2+n_1Dn_1+n_2Dn_1+n_1Dn_2+n_2Dn_2),\\[1mm]
f_3\sim &\D\mathcal{O}(1)(w_2Dn_1+w_1Dn_2+n_1Dw_2+n_2Dw_1),\\[1mm]
f_4\sim &\D
\mathcal{O}(1)(w_1Dw_2+w_2Dw_1+n_1Dn_1+n_2Dn_1+n_1Dn_2+n_2Dn_2).
\end{array}
\eqno(3.7)
$$

Thus, we have
$$
 \begin{array}[b]{rl}
&\D\|(f_1,f_2)\|_{L^1}(s)\\[1mm]
\leq\!\! &\D
\!\!C\|(w_2Dn_1,w_1Dn_2,n_1Dw_1,n_2Dw_2,w_1Dw_1,w_2Dw_2,n_1Dn_1,n_2Dn_1,n_1Dn_2,n_2Dn_2)\|_{L^1}\\[1mm]
\leq\!\! &\D
\!\!C\{\|w_2\|\|Dn_1\|_{L^2}\!+\!\|w_1\|\|Dn_2\|\!+\!\|n_1\|\|Dw_1\|\!+\!\|n_2\|\|Dw_2\|\!+\!\|n_1\|\|Dw_2\|\!+\!\|n_2\|\|w_1\|\\[1mm]
&\D\!+\!\|w_1\|\|Dw_1\|\!+\!\|w_2\|\|Dw_2\|\!+\!\|n_1\|\|Dn_1\|\!+\!\|n_2\|\|Dn_1\|\!+\!\|n_1\|\|Dn_2\|\!+\!\|n_2\|\|Dn_1\|\}\\[1mm]
\leq\!\! &\D\!\! CM^2(t)(1+s)^{-2}.
\end{array}
\eqno(3.8)
$$

In fact, we just need to consider the key term
$\|n_1Dn_1\|_{L^1}(s)$ with the ``worst" decay rate (see ansatz
(3.4)), which can be estimated as
$$
\|n_1Dn_1\|_{L^1}(s)\leq C\|n_1\|\|Dn_1\|\leq
CM^2(t)(1+s)^{-\frac{3}{4}}(1+s)^{-\frac{5}{4}}=CM^2(t)(1+s)^{-2}.\eqno(3.9)
$$

Then, for $\|(f_1,f_2)\|(s)$, we get
$$
 \begin{array}[b]{rl}
&\D\!\!\|(f_1,f_2)\|(s)\\[1mm]
\leq &\D\!\!
C\|(w_2Dn_1,w_1Dn_2,n_1Dw_1,n_2Dw_2,w_1Dw_1,w_2Dw_2,n_1Dn_1,n_2Dn_1,n_1Dn_2,n_2Dn_2)\|\\[1mm]
\leq &\D\!\!
C\{\|w_2\|_{\infty}\|Dn_1\|\!+\!\|w_1\|_{\infty}\|Dn_2\|\!+\!\|n_1\|_{\infty}\|Dw_1\|\!+\!\|n_2\|_{\infty}\|Dw_2\|\!+\!\|n_1\|_{\infty}\|Dw_2\|\\[1mm]
&\D\ \ \ \ \ \!\!+\!\|n_2\|_{\infty}\|w_1\|
\!+\!\|w_1\|_{\infty}\|Dw_1\|\!+\!\|w_2\|_{\infty}\|Dw_2\|\!+\!\|n_1\|_{\infty}\|Dn_1\|\!+\!\|n_2\|_{\infty}\|Dn_1\|\\[1mm]
&\D\ \ \ \ \
\!\!+\!\|n_1\|_{\infty}\|Dn_2\|\!+\!\|n_2\|_{\infty}\|Dn_1\|\}.
\end{array}
\eqno(3.10)
$$

The key term $\|n_1Dn_1\|$ in (3.10) with the ``worst" decay rate
(see ansatz (3.4)) can be estimated as
$$\begin{array}[b]{rl}
\|n_1Dn_1\|(s)\leq &\D C\|n_1\|_{\infty}\|Dn_1\|\leq
C\|Dn_1\|^{\frac{3}{2}}\|D^2n_1\|\\[1mm]
\leq&\D
CM^2(t)(1+s)^{-(\frac{3}{2}\cdot\frac{5}{4}+\frac{1}{2}\cdot\frac{5}{4})}=CM^2(t)(1+s)^{-\frac{5}{2}},
\end{array}
\eqno(3.11)
$$
where we have used the fact (the Nirenberg inequality)
$$
\|u\|_{\infty}\leq
C\|Du\|^{\frac{1}{2}}\|D^2u\|^{\frac{1}{2}}.\eqno(3.12)
$$

Thus from (3.6) to (3.11), we get
$$
\|n_1\|(t)\leq
\epsilon_1(1+t)^{-\frac{3}{4}}+CM^2(t)(1+t)^{-\frac{3}{4}},\eqno(3.13)
$$
and
$$
\|w_1\|(t)\leq
\epsilon_1(1+t)^{-\frac{5}{4}}+CM^2(t)(1+t)^{-\frac{5}{4}}.\eqno(3.14)
$$

Similarly, we can easily obtain
$$
\|(f_3,f_4)\|_{L^1}(s)\leq CM^2(t)(1+s)^{-2},\eqno(3.15)
$$
and
$$
\|(f_3,f_4)\|(s)\leq CM^2(t)(1+s)^{-\frac{5}{2}},\eqno(3.16)
$$
which together with (3.6) yield that
$$
\|(n_2,w_2)\|(t)\leq \epsilon_1e^{-bt}+CM^2(t)(1+t)^{-2},\
b>0.\eqno(3.17)
$$

\bigbreak \noindent\textbf{Step 2. Estimates for higher order
derivatives} \ Firstly, we consider the first order derivatives.
From (3.7), we get
$$
 \begin{array}[b]{rl}
&\D\!\!\|D(f_1,f_2)\|(s)\\[1mm]
\leq &\D\!\!
C\|D(w_2Dn_1,w_1Dn_2,n_1Dw_1,n_2Dw_2,w_1Dw_1,w_2Dw_2,n_1Dn_1,n_2Dn_1,n_1Dn_2,n_2Dn_2)\|\\[1mm]
\leq &\D\!\!
C\{\|(Dw_2Dn_1,w_2D^2n_1,Dw_1Dn_2,w_1D^2n_2,Dn_1Dw_1,n_1D^2w_1,Dn_2Dw_2,n_2D^2w_2,\\[1mm]
&\D
Dw_1Dw_1,w_1D^2w_1,Dw_2Dw_2,w_2D^2w_2,Dn_1Dn_1,n_1D^2n_1,\\[1mm]
&\D Dn_1Dn_2,n_2D^2n_1,n_1D^2n_2,Dn_2Dn_2,n_2D^2n_2)\|\}.
\end{array}
\eqno(3.18)
$$

In the same way, from ansatz (3.4), we only consider the terms
containing $n_1$. In fact, from (3.12) we have
$$
\begin{array}[b]{rl}
&\D \|n_1D^2n_1\|(s)\leq C\|n_1\|_\infty\|D^2n_1\|\leq
C\|Dn_1\|^{\frac{1}{2}}\|D^2n_1\|^{\frac{3}{2}}\\[1mm]
 \leq &\D
 CM^2(t)(1+s)^{-(\frac{1}{2}\cdot\frac{5}{4}+\frac{3}{2}\frac{5}{4})}=CM^2(t)(1+s)^{-\frac{5}{2}},
 \end{array}
 \eqno(3.19)
$$
$$
\begin{array}[b]{rl}
&\D \|Dn_1Dn_1\|(s)\leq C\|Dn_1\|_\infty\|Dn_1\|\leq
C\|D^2n_1\|^{\frac{1}{2}}\|D^3n_1\|^{\frac{1}{2}}\|Dn_1\|\\[1mm]
\leq &\D
CM^2(t)(1+s)^{-(\frac{1}{2}\cdot\frac{5}{4}+\frac{5}{4})}=CM^2(t)(1+s)^{-\frac{15}{8}}.
 \end{array}
 \eqno(3.20)
$$

Consequently, from (3.6) and (3.8), we have
$$
\|Dn_1\|(t)\leq
\epsilon_1(1+t)^{-\frac{5}{4}}+CM^2(t)(1+t)^{-\frac{15}{8}},\eqno(3.21)
$$
and
$$
\|Dw_1\|(t)\leq
\epsilon_1(1+t)^{-\frac{7}{4}}+CM^2(t)(1+t)^{-\frac{15}{8}}.\eqno(3.22)
$$

Similarly, (3.19) and (3.20) together with (3.6) and (3.8) yield
that
$$
\|D(n_2,w_2)\|(t)\leq
\epsilon_1e^{-bt}+CM^2(t)(1+t)^{-\frac{15}{8}},\ b>0.\eqno(3.23)
$$

Next, we can use the same method as above to derive the estimates of
 second order derivatives. from (3.7), we have
$$\begin{array}[b]{rl}
\|D^2(f_1,f_2)\|(s)= &\D\!\|(D^2w_2Dn_1,Dw_2D^2n_1,w_2D^3n_1,D^2n_1Dw_1,Dn_1D^2w_1,\\[1mm]
 &\D \! D^2n_1Dn_1,n_1D^3n_1,D^2n_1Dn_2,Dn_1D^2n_2,n_1D^3n_2,\cdots)\|(s),
 \end{array}
 \eqno(3.24)
$$
where ``$\cdots$" denotes the other terms containing no variable
$n_1$. Indeed, because of the ``worst" decay rate of $n_1$ in the
ansatz (3.4), we just need to consider the terms listed above,
especially two terms $\|Dn_1D^2n_1\|$ and $\|n_1D^3n_1\|$. In fact,
we have
$$
\|Dn_1D^2n_1\|(s)\leq C\|Dn_1\|_\infty\|D^2n_1\|\leq
C\|D^2n_1\|^{\frac{3}{2}}\|D^3n_1\|^{\frac{1}{2}}\leq
CM^2(t)(1+s)^{-\frac{15}{8}},\eqno(3.25)
$$
and
$$\begin{array}[b]{rl}
\|n_1D^3n_1\|(s)\leq &\D C\|n_1\|_\infty\|D^3n_1\|\leq
C\|Dn_1\|^{\frac{1}{2}}\|D^2n_1\|^{\frac{1}{2}}\|D^3n_1\|\\[1mm]
\leq &\D
CM^2(t)(1+s)^{-(\frac{1}{2}\cdot\frac{5}{4}+\frac{1}{2}\cdot\frac{5}{4})}=CM^2(t)(1+s)^{-\frac{5}{4}}.
\end{array}
\eqno(3.26)
$$

Then, (3.25), (3.26), (3.6) and (3.8) imply that
$$
\|D^2(n_1,w_1)\|(t)\leq
\epsilon_1(1+t)^{-\frac{5}{4}}+CM^2(t)(1+t)^{-\frac{5}{4}}.\eqno(3.27)
$$

Similarly, we have
$$
\|D^2(n_2,w_2)\|(t)\leq
\epsilon_1(1+t)^{-\frac{9}{4}}+CM^2(t)(1+t)^{-\frac{5}{4}}.\eqno(3.28)
$$

Lastly, Theorem 1.1 and the relation (2.2) lead to that
$$
\|D^3(n_1,n_2,w_1,w_2)\|(t)\leq C\epsilon_0.\eqno(3.29)
$$

In summary, from (3.13), (3.14), (3.17), (3.21)-(3.23) and
(3.27)-(3.29), we deduce that
$$
M(t)\leq C(\epsilon_0+\epsilon_1)+CM^2(t).\eqno(3.30)
$$

By the standard continuous argument, we know that there exists a
constant $C>0$ such that
$$
M(t)\leq C(\epsilon_0+\epsilon_1), \ t\in[0,T].\eqno(3.31)
$$

This proves Theorem 1.2. Immediately, Corollary 1.1 can be derived
from Theorem 1.2 and (2.2).
\\
\\
{\bf Acknowledgement:} \ \   The research of Z.G. Wu was supported
by NSFC (No. 11101112 and 11071162). The research of Y.M. Qin was
supported by NSFC (No. 11271066, No. 11031003) and Innovation
Program of Shanghai Municipal Education Commission ( No. 13ZZ048).

\bibliographystyle{plain}

\end{document}